\documentclass[11pt]{article}
\usepackage{amsfonts}

\usepackage{graphics}
\usepackage{indentfirst}
\usepackage{cite}
\usepackage{latexsym}
\usepackage{amsmath}
\usepackage{amssymb}
\usepackage[dvips]{epsfig}
\usepackage{amscd}

\newtheorem{theorem}{Theorem}[section]
\newtheorem{remark}{Remark}[section]

\newtheorem{definition}{Definition}[section]
\newtheorem{lemma}[theorem]{Lemma}

\newtheorem{corollary}[theorem]{Corollary}
\newtheorem{proposition}[theorem]{Proposition}

\newcommand{\n}{\rho}

\newcommand{\ti}{\tilde}

\newcommand{\lm}{\lambda}

\def\pf{{\it Proof.}  }
\renewcommand{\div}{ {\rm div }  }

\newcommand{\na}{\nabla }

\newcommand{\pa}{\partial}

\newcommand{\bi}{\bibitem}

\newcommand{\bt}{\begin{theorem}}
\newcommand{\bl}{\begin{lemma}}
\newcommand{\el}{\end{lemma}}
\newcommand{\et}{\end{theorem}}
\newcommand{\ga}{\gamma}

\newcommand{\de}{\delta}
\newcommand{\ve}{\varepsilon}
\newcommand{\la}{\label}

\newcommand{\bn}{\begin{eqnarray}}
\newcommand{\en}{\end{eqnarray}}
\newcommand{\bnn}{\begin{eqnarray*}}
\newcommand{\enn}{\end{eqnarray*}}

\newcommand{\bnnn}{\begin{eqnarray*}}
\newcommand{\ennn}{\end{eqnarray*}}
\newcommand{\ben}{\begin{enumerate}}
\newcommand{\een}{\end{enumerate}}

\newcommand{\ba}{\begin{aligned}}
\newcommand{\ea}{\end{aligned}}
\newcommand{\be}{\begin{equation}}
\newcommand{\ee}{\end{equation}}

\def\p{\partial}
\def\norm[#1]#2{\|#2\|_{#1}}

\def\lap{\triangle}

\def\lam{\lambda}

\def\rrr{\mathbb{R}^3}
\def\xl{\left}
\def\xr{\right}

\def\O{\Omega}

\def\al{\alpha}
\def\tn{\n_\infty}
 
\makeatletter      
\@addtoreset{equation}{section}
\makeatother       

\title{On  local strong and  classical solutions   to the three-dimensional barotropic compressible Navier-Stokes equations with vacuum}
\date{}
\author{Xiangdi  H{\small UANG} \thanks{Institute of  Mathematics, AMSS,
Chinese Academy of Sciences, Beijing 100190, People's  Republic of
China ({\tt xdhuang@amss.ac.cn}).X.-D. Huang is partially supported by National Natural Science Foundation of China, Grant Nos. 11688101, 11731007,11671412 and Youth Innovation Promotion Association CAS.}
  }

\begin{document}
\maketitle

\begin{abstract} We consider  the local well-posedness of strong and classical solutions to the three-dimensional barotropic compressible Navier-Stokes equations     with  density containing vacuum initially. We first prove the local existence and uniqueness of the strong solutions, where the initial compatibility condition proposed in \cite{cho1,K2,coi1}
 is removed under suitable sense. Then, the continuous of strong solutions on the initial data is derived under an additional compatibility condition.  Moreover, for the initial data satisfying   some additional regularity  and   compatibility   condition, the strong solution is proved to be  a classical one.
\end{abstract}


\section{Introduction and main results}
  We consider the three-dimensional
  barotropic compressible Navier-Stokes equations    which read as follows:
\be\la{n1}
\begin{cases} \rho_t + \div(\rho u) = 0,\\
 (\rho u)_t + \div(\rho u\otimes u) + \nabla P  = \mu\lap u + (\mu + \lam)\nabla \div u ,
\end{cases}
\ee
where   $t\ge 0, x=(x_1,x_2,x_3)\in  \Omega\subset \rrr, \rho=\n(x,t),$ $u=(u_1(x,t),u_2(x,t),u_3(x,t)),$ and $P=P(\n),$  represent, respectively,  the density, the velocity, and the pressure. The constant viscosity  coefficients  $\mu$ and  $\lambda$ satisfy  the  physical hypothesis:
\be\la{n3}
\mu>0,\quad 2\mu+3\lam\ge 0.
\ee

 Let $\Omega\subset \rrr$ be either a smooth bounded domain or the whole space $\rrr$,  we impose the following initial and boundary conditions on   \eqref{n1}:
\be \la{n4} \n(x,0)=\n_0(x), \quad \n u(x,0)=m_0(x),\quad x\in \Omega,\ee
and
\be \la{n5}\begin{cases} u(x,t)=0,~~x\in\p \O,~~&\mbox{if}~\O  \subset \subset \rrr ,\\
(\n, u)(x,t)\rightarrow (\n_\infty, 0),~~\mbox{as}~|x|\rightarrow\infty,~~&\mbox{if}~\O=\rrr,
\end{cases}\ee
with constant $\n_\infty \ge 0.$

It is important to investigate the well-posedness of strong solutions for compressible Navier-Stokes equations. 

As long as the initial density is away from vacuum, the local well-posedness theory to the problem \eqref{n1} are established in \cite{vlli} and \cite{Na,se1}, respectively.  In 1980s, Matsumura-Nishida \cite{M1} proved the existence of global classical solutions when the initial data are 
close to a non-vacuum resting states. Besides, it is shown by Hoff \cite{Ho4,Hof2} that the system will admit at least one global weak solution with strictly positive initial density and temperature for discontinuous initial data.
 
 Things become more complicated when the density is allowed to vanish. In 1994, The major breakthrough is due to Lions \cite{L4,L1} (then improved by  Feireisl \cite{Fe,F1}), where global existence of weak solutions with finite energy without any size restriction on the initial data can be proved under the condition that the exponent $\ga$ is suitably large.
 Later, Hoff \cite{Ho3,hs,ht} obtained a new type of global weak solutions with small energy. 
 Considering the strong or classical solutions with vacuum, the authors in \cite{cho1,K2,coi1, sal}  obtained the local existence and uniqueness of strong and classical solutions  for  three-dimensional bounded or unbounded domains and for two-dimensional bounded ones. It should be noted that the results in those of \cite{cho1,K2,coi1, sal} are derived under some additional  compatibility conditions, see \eqref{co2} in the below.  More precisely, they required that $g\in L^2(\O)$ or $g\in H^1(\O)$ in \eqref{co2} for the strong or classical solutions, respectively.  In this direction, a natural question arises whether one can remove or relax the initial compatibility conditions with nonnegative density in suitable sense. Indeed, this is the aim of this paper, i.e, we establish the local existence of strong solutions without the initial compatibility condition.

Before stating the main results, we first explain the notations and
conventions used throughout this paper. For $1\le r\le \infty $ and $k\ge 1$,    the standard Lebesgue and
  Sobolev spaces are defined as follows:
   \bnn  \begin{cases}L^r=L^r(\Omega ),\quad
W^{k,r}  = W^{k,r}(\Omega) , \quad H^k = W^{k,2} ,\\ D_0^1=\begin{cases}H_0^1(\O),&\mbox{ for   bounded }\O\subset \rrr,\\ \{f\in L^6|\na f\in L^2\}&\mbox{ for   }\O= \rrr.\end{cases}\end{cases}\enn  The first main result of this paper  is  the  following Theorem \ref{t1} concerning the  local existence   of  strong solutions whose definition is as follows:
\begin{definition} If  all derivatives involved in \eqref{n1} for $(\rho,u)  $  are regular distributions, and   equations  \eqref{n1} hold   almost everywhere   in $\O\times (0,T),$ then $(\n,u)$  is called a  strong solution to  \eqref{n1}.
\end{definition}

\begin{theorem}\la{t1} Assume that $P=P(\cdot)\in C^1[0,\infty).$ For some  $3<q<6 $ and $\n_\infty\ge 0$,  assume that  the initial data $( \n_0 ,m_0)$ satisfy
  \be\la{1.9}\n_0\ge 0,\,\,
 \rho_0-\n_\infty\in    L^{\ti p}\cap D^1\cap W^{1,q},\,\, u_0 \in D_0^1,
   \ee and \be \la{a1.9} m_0=\n_0u_0,\ee  where \be \ti p\triangleq\begin{cases}3/2, &\mbox{ for } \O=\rrr \mbox{ and } \tn=0,\\  2, &\mbox{ otherwise. }  \end{cases}\ee    Then there exists a positive time $T_0>0$ such that the problem  \eqref{n1}--\eqref{n5} has a unique   strong solution $(\n,u)$ on $\O\times (0,T_0]$ satisfying that
  \be\la{1.10}\begin{cases}
  \rho-\n_\infty\in C([0,T_0];  L^{\ti p}\cap D^1 \cap W^{1,q} ),    \\  \na u,\,  \sqrt{t}\na^2 u,
   \,    \sqrt{t} \sqrt{\n}  u_t,\, t \na u_t \in L^\infty(0,T_0; L^2 ) ,   \\
  {t}\na u\in L^\infty(0,T_0; W^{1,q} )   ,  \, \sqrt{\n} u_t, \, \sqrt{t}\na u_t  \in L^2(\O\times(0,T_0)).
   \end{cases}\ee  Furthermore, if in addition to \eqref{1.9} and \eqref{a1.9},  $(\n_0,u_0)$ satisfies the  compatibility conditions
\be \la{co2}- \mu\lap u_{0 } - (\mu + \lam)\nabla \div u_{0 } +  \nabla P(\n_{0 })=\n_0^{1/2}g ,     \ee for  $g \in L^2  ,$   $(\n,u)$ also satisfies \be \la{lco1} \begin{cases}\na u\in L^\infty(0,T_0;H^1),\, \sqrt{t}\na u\in L^\infty(0,T_0;W^{1,q}),\\ \sqrt{\n} u_t, \sqrt{t}\na u_t  \in L^\infty(0,T_0;L^2),\,\na u_t\in L^2(\O\times (0,T_0))    .\end{cases}\ee
\end{theorem}

Next, the following Corollary \ref{t2} whose proof is similar as that of \cite[Theorem 3]{K2} gives  the continuous
dependence of the solution on the data   provided  \eqref{co2} holds.

\begin{corollary}\la{t2} For each $i=1,2,$ let $(\n_i,u_i)  $ be the local strong solution to the   problem  \eqref{n1}--\eqref{n5}  with the initial data $(\n_{0i},u_{0i})$   satisfying \eqref{1.9}, \eqref{a1.9},  and the  compatibility conditions \eqref{co2} with $g=g_i.$ Moreover, assume that  $(\n_{0i}, u_{0i})$  satisfies
\be\label{ma1} \|\n_{0i}-\tn\|_{ L^{\ti p} \cap D^1\cap W^{1,q}}+\|\na u_{0i}\|_{H^1}+ \|g_{ i}\|_{L^2}\le K.\ee
Then  there exists a small time $T_0$ and a positive constant $C$ depending only on $T_0$ and $ K$
such that
\be\ba\label{ma2}
&\sup\limits_{0\le t\le T_0} \xl(\|\rho_1^{1/2}(u_1-u_2) \|_{L^2}^2   +\|\n_{1}-\n_{2}\|_{L^{\ti p}}^2\xr)
 +  \int_0^{T_0} \|\na (u_1-u_2)\|_{L^2}^2ds \\
 &\le C\|\n_{01}^{1/2}(u_{01}-u_{02})\|_{L^2}^2+C\|\n_{01}-\n_{02}\|_{L^{\ti p}}^2
  .
\ea\ee
    \end{corollary}

Finally, if the initial data $(\n_0,m_0)$ satisfy some additional regularity and compatibility conditions, the local strong solution $(\n,u) $ obtained by Theorem \ref{t1} becomes a classical one.
\begin{theorem}\la{t3} Assume that $P(\n)$ satisfies either \be P(\cdot)\in C^2[0,\infty)\ee or \be P(\n)=A\n^\ga (A>0,\ga>1).\ee In addition to    \eqref{1.9}, \eqref{a1.9}, and \eqref{co2},  assume further that
  \be\la{1.c1}  \na^2 \n_0,\, \na^2 P(\n_0 )\in L^2\cap  L^q .
 \ee    Then,   in addition to \eqref{1.10} and \eqref{lco1},    the  strong  solution  $(\rho,u)$ obtained by Theorem \ref{t1} satisfies    \be\la{1.a10}\begin{cases}
  \na^2\rho,   \,\,\na^2 P(\rho)\in C([0,T_0];L^2\cap L^q  ), \\       \na u\in L^2(0,T_0;H^2),\,
  \sqrt{t}\na u\in L^\infty(0,T_0;H^2),\\ t\na u\in L^\infty(0,T_0;W^{2,q}), \, \sqrt{t}\na  u_t\in L^2(0,T_0;H^1),\\t \na u_t\in L^\infty(0,T_0;H^1),\quad
 t u_{tt}\in L^2 (0,T_0;D_0^1),\\
  t\sqrt{\n}u_{tt}\in L^\infty(0,T_0;L^2),\quad \sqrt{t}\sqrt{\n}u_{tt}\in L^2(0,T_0;L^2).\end{cases}\ee
      \end{theorem}

A few remarks are in order:
\begin{remark}  To obtain the local existence and uniqueness of strong solutions, in Theorem \ref{t1}, the only compatibility condition we   need is \eqref{a1.9} which is   much weaker than those of   \cite{cho1,
K2,coi1, sal}  where   not only \eqref{a1.9} but also \eqref{co2}  is   needed.  Moreover, the strong solutions obtained in Theorem \ref{t1} are somewhat more regular   than those in \cite{cho1,K2,coi1} when $t>0$. In this sense, we successfully remove the  compatibility condition required in   \cite{cho1,
K2,coi1, sal}.
 \end{remark}

\begin{remark}  After obtaining  the existence result in Theorem \ref{t1}, the  continuous
dependence of the solution on the data is shown in Corollary  \ref{t2}, provided that  the initial data satisfy the compatibility condition  \eqref{co2}. Indeed, Theorem  \ref{t1} and Corollary   \ref{t2} tell us how the  the compatibility condition  \eqref{co2} plays its role in discussing the local well posedness of strong solutions to the problem (\ref{n1})-(\ref{n5})  with vacuum.
 \end{remark}

 \begin{remark}  For the local existence of classical solutions obtained in Theorem \ref{t3}, we only need the initial data satisfying  the compatibility condition  \eqref{co2} for some $g\in L^2$ which is  in sharp contrast to  Cho-Kim \cite{coi1} where  the compatibility condition  \eqref{co2}  is needed for $g\in H^1$.  This means that our Theorem \ref{t3} essentially weaken those assumptions on the compatibility condition in \cite{coi1}.
 \end{remark}

We now comment  on the analysis of this paper.
 First, we will consider the approximating system for the initial density strictly away from vacuum,
 whose local existence theory has been shown in Lemma \ref{th0}.
 By employing some basic ideas due to Hoff \cite{Ho4,Hof2} and careful analysis,
 we succeed in deriving the uniform a priori estimates on the density and velocity
 which are independent of the lower bound of the density.
To do this,  the key issue is to get the uniform upper bound of the density without requiring the additional compatibility condition  \eqref{co2}.  Indeed, this is achieved by deriving the time weighted estimates on $\|\sqrt{\n}u_t\|_{L^2}$ and $\|\na u_t\|_{L^2}$, see Lemma \ref{l3.2},  which are crucial for bounding the $L^1L^\infty$-norm of $\na u$ and thus
getting the uniform upper bound of the density.
Then, with the desired estimates on solutions at hand, we will apply the standard compact arguments which show that the limit is exactly the strong solutions of the original one.
Finally, for the initial data satisfying some additional regularity and compatibility conditions, the standard arguments will be used to obtain
the higher order estimates of the
solutions which are needed to guarantee the local strong solution to be
a classical one.

  We shall briefly describe the structure of this article. Some fundamental Lemmas will be exhibited in section 2. To get the local existence and uniqueness of strong and classical solutions, some a priori estimates in section 3 and 4 are established in orders. Consequently, we arrive the results of Theorems \ref{t1} and \ref{t3} in Section 5.

\section{Preliminaries}\la{sec2}

First, in   this section and the following two,  we denote
 \be \la{or1} \O_R=\begin{cases}\O, &\mbox{\rm for bounded }\,\O\subset\rrr,\\ B_R\triangleq\{x\in \rrr||x|<R\}, &\mbox{\rm for   }\,\O=\rrr,
 \\ \end{cases}\ee
and
   \bnn L^p=L^p(\O_R ) , \quad
W^{k,p}  = W^{k,p}(\O_R) , \quad H^k = W^{k,2}, \enn for $p\ge 1$ and positive integer $k.$

Then, for the initial
density    strictly away from vacuum,   the following    local existence theory can be shown by similar arguments as in   \cite{cho1,K2,coi1,vlli}.
\begin{lemma}\la{th0}  Assume that $P(\cdot)\in C^3[0,\infty)$ and that
 the initial data $(\n_0 ,m_0)$ satisfy
 \bnn 0<\delta\le \n_0,~~
 \rho_0 \in   H^3,\,\,
 u_0 \in H_0^1\cap H^3,\,\, m_0=\rho_0 u_0 .\enn
   Then
there exist  a small time $T_*>0$  such that   the problem  \eqref{n1}--\eqref{n5} admits a unique classical solution
$(\rho, u)$
on $\O_R\times(0,T_*]$ satisfying
\bnn\begin{cases}    \n\in C\left([0,T_*];H^{3}\right), u\in C\left([0,T_*]; H^1_0\cap H^{3}\right)\cap L^{2}\left(0,T_*;H^{4}\right),\\
 u_{t}\in L^{\infty}\left(0,T_*;H_0^1\right)\cap  L^{2}\left(0,T_*;H^{2}\right),\sqrt{\n}u_{tt}\in L^{2}\left(0,T_*;L^{2}\right),\\
 \sqrt{t}  u \in L^{\infty}\left(0,T_*;H^{4}\right), \, \sqrt{t} u_{t}\in  L^{\infty}\left(0,T_*; H^{2}\right),\,\sqrt{t} u_{tt}\in L^{2}\left(0,T_*;H^{1}\right),\\
 \sqrt{t}\sqrt{\n}u_{tt}\in L^{\infty}\left(0,T_*;L^{2}\right), tu_{t}\in L^{\infty}\left(0,T_*;H^{3}\right), \\
  tu_{tt}\in L^{\infty}\left(0,T_*;H^{1}\right) \cap  L^{2}\left(0,T_*;H^{2}\right), t\sqrt{\n}u_{ttt}\in L^{2}\left(0,T_*;L^{2}\right),\\
 t^{3/2} u_{tt}\in  L^{\infty}\left(0,T_*;H^{2}\right),t^{3/2} u_{ttt}\in L^{2}\left(0,T_*;H^{1}\right), \\ t^{3/2}\sqrt{\n}u_{ttt}\in L^{\infty}\left(0,T_*;L^{2}\right).
 \end{cases}\enn
\end{lemma}

Next, the following well-known Gagliardo-Nirenberg inequality
  will be used later frequently (see \cite{la}).

\begin{lemma}
[Gagliardo-Nirenberg]\la{l1} For  $p\in [2,6],q\in(1,\infty), $ and
$ r\in  (3,\infty),$ there exists some generic
 constant
$C>0$ independent of $R$ such that for   $f\in H^1_0(\O_R) $
and $g\in L^q (\O_R)\cap W^{1,r} (\O_R), $    \bn
\la{g1}\|f\|_{L^p}^p\le C \|f\|_{L^2}^{(6-p)/2}\|\na
f\|_{L^2}^{(3p-6)/2} ,\en  \bn
\la{g2}\|g\|_{L^\infty} \le C
\|g\|_{L^q}+ C
\|g\|_{L^q}^{q(r-3)/(3r+q(r-3))}\|\na g\|_{L^r}^{3r/(3r+q(r-3))} .
\en
\end{lemma}

Finally, we state the following $L^p$-bounds for the weak solutions to    the Lam\'e system with the Dirichlet boundary conditions
\be \la{lame}\ba\begin{cases}
-\mu\Delta v-(\mu+\lambda)\na \div v=F,\,\, &x\in\O_R,\\
 v=0, \,\, &x\in\pa\O_R.\end{cases}\ea\ee
\begin{lemma}[\cite{cho1,adn}]\label{ADN} For $p>1$ and $k\ge 0,$ there exists a positive constant $C$ independent of $R$ such that
\be \la{lp} \|\na^{k+2}v\|_{L^p(\O_R) }\le C \|F\|_{W^{k,p} (\O_R)} ,\ee
for every  solution  $v\in W_0^{1,p}(\O_R)$ of \eqref{lame}.
\end{lemma}

\section{A priori estimates (I)}


Let $\O_R$ and $(\n_0,m_0)$  be as in Lemma \ref{th0} and $(\n,u)$  the solution to  the   problem  \eqref{n1}--\eqref{n5} on $\O_R\times (0,T_*]$ obtained by Lemma \ref{th0}.  For  $q\in (3,6),$  we denote  \be \psi(t)\triangleq 1+ \|\na u\|_{L^2} +    \|  \n-\tn  \|_{ L^{\ti p}\cap D^{1}\cap W^{1,q}} .  \ee 
  Then the main aim of this section is to derive
  the following key a priori estimate  on $\psi . $
\begin{proposition} \la{pro}   For $q\in (3,6),$  there exist positive constants $T_0$ and $M$ both  depending only on   $\mu$, $\lam$, $P,$ $q,$ $\tn,$ $\psi(0),$  and $\O$ but independent of   $ R$  such that  \be\la{o1}\ba &\sup\limits_{0\le t\le T_0}\xl(\psi(t)+t(\|\na^2u\|_{L^2}^{2}+\|\sqrt{\n} u_t\|_{L^2}^{2})+t^2(\|\na u_t\|_{L^2}^{2}+\|\na^2u\|_{L^q}^{2})\xr) \\&+\int_0^{T_0}t\|\na u_t\|_{L^2}^2dt \le M.\ea\ee
 \end{proposition}

To prove Proposition \ref{pro},  we begin with the following      $L^2$-bound  for $\nabla u.$

\begin{lemma} \la{l3.0}     There exist positive constants $\alpha=\alpha(  q)>1$  such that
\be\ba\la{3.1}
  &\sup_{0\le s\le t} (\|\na u\|^2_{L^{2}}+\|P-P(\tn)\|^2_{L^2})
  +\int_{0}^{t}\|\sqrt{\rho} u_t\|_{L^2}^2  ds  \\& \le C  +C\int_{0}^{t} M_P(\psi) \psi^{\alpha} ds,
\ea
\ee  where and in this section,
 \be M_P(\psi)\triangleq 1+\max_{0\le s\le \psi}(|P (s)|+|P'(s)|) ,\ee and $C$ denotes a
generic positive constant
 depending only on       $\mu$, $\lam$, $P, $ $q,$   $\tn,$   $ \psi(0),$   and $\O$ but independent of $R.$
\end{lemma}

\pf   First,  multiplying    equations  $\eqref{n1}_{2}$  by  $ u_t$ and    integrating the resulting equations by parts yield
\be\la{3r1}\ba& \frac{d}{dt}\int\left((\mu+\lm)(\div u)^2+\mu|\na u|^2\right)dx+\int \n|  u_t|^2dx\\ &\le C\int\n|u|^2|\na u|^2dx +2\int (P-P(\tn))\div u_tdx ,\ea\ee where, in this section and the next, we denote $$\int \cdot dx=\int_{\O_R}\cdot dx.$$

Then, on the one hand,   the Gagliardo-Nirenberg inequality implies  that
\be\ba\la{cc5}\int\n |u|^2|\na u|^2dx &\le  \|\n \|_{L^\infty} \|u\|_{L^6}^2\|\na u\|_{L^3}^{2} \\&\le C\|\n \|_{L^\infty} \|\na u\|_{L^2}^3\|\na  u\|_{H^1} \\&\le C  \psi^{\alpha }\|\na^2 u\|_{L^2}+ C  \psi^{\alpha},\ea\ee
where (and in what follows)   $\alpha=\alpha(  q)>1.$
Note that $u $ is a solution of the following elliptic system
\be\label{ell1}\ba \begin{cases}
- \mu\lap u - (\mu + \lam) \nabla \div u =-\n  ( u_t+u\cdot\na u) -\nabla P ,   \quad  & x\in \O_R,  \\
 u  =0,  \quad   &x\in\pa\O_R.  \end{cases}\ea\ee
 Applying Lemma \ref{ADN} to \eqref{ell1} yields
 \bnn   \ba  \|\na^2u\|_{L^2} \le & C \left(\|\n(   u_t  + u\cdot\na u)\|_{L^2}  + \|\na P\|_{L^2}\right)  \\
\le & C\psi^{1/2}\|\sqrt{\n} u_t\|_{L^2}+CM_P(\psi)\psi^\alpha +\frac{1}{2} \|\na^2u\|_{L^2},\ea\enn
where in the second inequality we have used   \eqref{cc5}. This implies \be \label{ell2} \|\na^2u\|_{L^2} +\|\n(   u_t  + u\cdot\na u)\|_{L^2}
\le   C\psi^{1/2}\|\sqrt{\n} u_t\|_{L^2}+CM_P(\psi)\psi^\alpha  .\ee
On the other hand,
we deduce from   the Sobolev inequality that
\be\la{bv3}\ba &2\int (P-P(\tn))\div u_tdx\\
&=2\frac{d}{dt}\int (P-P(\tn))\div u dx-2\int  P'(\n)\n_t\div udx   \\
&\le 2\frac{d}{dt}\int (P-P(\tn))\div udx+C M_P(\psi)\psi^2,\ea\ee where we have used
\be\la{p01}\|\n_t\|_{L^2}\le C\|u\|_{L^6}\|\na \n\|_{L^3}+C\|\n\|_{L^\infty}\|\na u\|_{L^2}\le C\psi^2,\ee due to \eqref{n1}$_1.$

Substituting \eqref{cc5},   \eqref{ell2}, and  \eqref{bv3}  into \eqref{3r1} and using Cauchy's inequality lead to
   \be\label{ma3.1} \ba   &\frac{d}{dt}\int\left((\mu+\lm)(\div u)^2+\mu|\na u|^2- 2(P-P(\tn))\div u\right)dx+\int \n|  u_t|^2dx\\ &\le C \psi^{\alpha}\|\n^{1/2}u_t\|_{L^2}+CM_P(\psi) \psi^\alpha\\
   &\le \frac{1}{2} \|\n^{1/2}u_t\|_{L^2}^2+C  M_P(\psi) \psi^\alpha.\ea\ee

   Finally, it follows from \eqref{p01} that\be \la{p3.15}\ba \frac{d}{dt}\|P-P(\tn)\|^2_{L^2}&\le C\int|P-P(\tn)||P'(\n)||\n_t| dx\\&\le CM_P(\psi)\psi^\al,\ea\ee which together with \eqref{ma3.1} gives \eqref{3.1} and finishes
    the proof of  Lemma \ref{l3.0}.  \hfill $\Box$

\begin{lemma}\la{l3.2}It holds that
 \be\ba\la{li-1a}
 \sup_{0\le s\le t} s\int \n |u_{t}|^2dx+\int_{0}^{t}s\|\na u_t\|_{L^2}^2 ds \le C   \exp\left\{C \int_0^tM_P^2(\psi)\psi^\alpha ds\right\}.\ea \ee
\end{lemma}

{\it Proof.}  Differentiating $\eqref{n1}_2$ with respect to $t$ gives
\be\la{zb1}\ba &-\mu\Delta u_t-(\mu+\lm)\na \div u_t\\&=-\n u_{tt}-\n u\cdot \na u_t-\n_t(u_t+u\cdot\na u)-\n u_t\cdot\na u-\na P_t.\ea\ee
 Multiplying \eqref{zb1} by $u_t,$  we obtain after using  integration by parts and $\eqref{n1}_1$ that\be\ba  \la{na8}&\frac{1}{2}\frac{d}{dt} \int \n |u_t|^2dx+\int \left((\mu+\lm)(\div u_t)^2+\mu |\na u_t|^2\right)dx\\
 &=-2\int \n u \cdot \na  u_t\cdot u_tdx  -\int \n u \cdot\na (u\cdot\na u\cdot u_t)dx\\
  &\quad-\int \n u_t \cdot\na u \cdot  u_tdx +\int P_{t}\div u_{t} dx\\
 &\le C\int  \n |u||u_{t}| \left(|\na  u_t|+|\na u|^{2}+|u||\na^{2}u|\right)dx +C\int \n |u|^{2}|\na u ||\na u_{t}|dx \\
 &\quad+C\int \n |u_t|^{2}|\na u |dx  +C\int |P_{t}||\div u_{t}|dx \triangleq \sum_i^4 J_i.
  \ea\ee

  We estimate each term  on the right-hand side of  \eqref{na8} as follows:

First,
it follows from    the Holder
  and the Gagliardo-Nirenberg inequalities    that
 \be\la{na2}\ba  
J_1 & \le C \|\n\|_{L^\infty}^{1/2} \| u\|_{L^{6}}\|\n^{1/2} u_{t}\|_{L^{2}}^{1/2} \|\n^{1/2} u_{t}\|_{L^{6}}^{1/2}\| \na u_{t}\|_{L^{2}}  \\
 &\quad +C\|\n\|_{L^\infty} \|u \|_{L^{6}} \|u_{t}\|_{L^{6}} \| \na u\|_{L^{3}}^{2}  +C\|\n\|_{L^\infty} \|u \|_{L^{6}}^{2}\|u_{t}\|_{L^{6}}  \| \na^{2} u \|_{L^{2}}  \\
 & \le C\psi^{\alpha} \|\n^{1/2} u_{t}\|_{L^{2}}^{1/2}  \| \na u_{t}\|_{L^{2}}^{3/2}+ C\psi^{\alpha}  \| \na u_{t}\|_{L^{2}}\| \na  u \|_{H^1}\\
 &\le  \ve\| \na u_{t}\|_{L^{2}}^{2}+C(\ve)\psi^{\alpha} \left(1+\| \na^{2} u \|_{L^{2}}^{2} +  \|\n^{1/2} u_{t}\|_{L^{2}}^{2} \right),
 \ea\ee
and
\be \ba \la{5.ap3}  J_2+J_3&\le C\|\n\|_{L^\infty} \|u\|_{L^{6}}^{2}\|\na u\|_{L^{6}} \| \na u_{t}\|_{L^{2}}+ C\| \na u\|_{L^{2}}
   \|\sqrt{\n}u_{t}\|_{L^{6}}^{3/2}\|\sqrt{\n}  u_{t}\|_{L^{2}}^{1/2} \\
&\le  \ve \| \na u_{t}\|_{L^{2}}^{2}+  C(\ve) \psi^{\alpha} \| \na^{2} u \|_{L^{2}}^{2}+ C(\ve)\psi^{\alpha}  \|\n^{1/2} u_{t}\|_{L^{2}}^{2}.\ea\ee

Next, it follows from  \eqref{p01}   that
 \be \la{na1}\ba  J_4
 &\le C   \|P'(\n)\|_{L^\infty}\|\n_t\|_{L^2}\|\na u_{t}\|_{L^{2}} \\  &\le \ve \|\na u_{t}\|_{L^{2}}^{2} + C(\ve)M_P^2(\psi)\psi^{\alpha} .\ea\ee

 Substituting \eqref{na2}--\eqref{na1} into \eqref{na8} and choosing $\ve$ suitably  small lead to
\be\ba\la{a4.6} & \frac{d}{dt} \int \n |u_t|^2dx+\int \left((\mu+\lm)(\div u_t)^2+\mu |\na u_t|^2\right)dx\\
  &  \le   C\psi^{\alpha}\left(1+\| \n^{1/2}u_{t}\|_{L^{2}}^{2}+\| \na^{2}u\|_{L^{2}}^{2}\right)\\
 &  \le   C \psi^{\alpha}  \|\n^{1/2} u_{t}\|_{L^{2}}^{2}+ C M_P^2(\psi)\psi^{\alpha},
  \ea\ee
 where  in the last inequality one has used \eqref{ell2}.

Finally, multiplying  \eqref{a4.6}  by $t , $  we obtain \eqref{li-1a} after
using Gronwall's  inequality and \eqref{3.1}. The proof of  Lemma \ref{l3.2} is completed.  \hfill $\Box$

\begin{lemma}\la{l3.4}It holds that
  \be\la{b1a} \ba  \sup\limits_{0\le s\le t}  \|\n-\tn \|_{L^{\ti p}\cap D^1 \cap  W^{1,q}}     \le  C\exp\left\{C\int_{0}^{t}M_P^2(\psi) \psi^{\alpha} ds \right\} . \ea \ee
\end{lemma}

 \emph{Proof}. First, using \eqref{n1}$_1,$   we have \be \la{p3.24}\frac{d}{dt}\|\n-\tn\|_{L^{\ti p}}\le C\psi^\al.\ee

 Next, differentiating \eqref{n1}$_1$ with respect to $x_i$ and multiplying the resulting equation by $r|\p_i \n|^{r-2}\na\n$ with $r\in [ 2,q]$, we obtain after  integration by parts that
 \be\label{ma3.2} \ba \frac{d}{dt}\|\na \n\|_{L^r} \leq  &
C\left(\|\nabla u \|_{L^\infty}\|\na \n\|_{L^r} +\|\n \|_{L^\infty}\|\na^2 u\|_{L^r}\right)\\ \leq  &
C\psi \left(\|\nabla u \|_{L^\infty}  + \|\na^2 u\|_{L^r}\right).
\ea \ee
Taking $r=2,q$ in \eqref{ma3.2} and using the Gagliardo-Nirenberg inequality, we have
 \bnn \ba \frac{d}{dt}\|\na \n\|_{L^2\cap L^q}
\le  C (1+\|\na^2 u\|_{L^2\cap L^q})\psi^\al ,
\ea \enn
 which together with  \eqref{p3.24}    yields \eqref{b1a} provided we show that
  \be\ba  \la{bb6-1}  \int_{0}^{t} \|\na^2u\|_{L^2\cap L^q}^{p_0}  ds
  \le  C   \exp\left\{C \int_0^tM_P^2(\psi)\psi^\alpha ds\right\}, \ea\ee for   \bnn \label{pp} p_0\triangleq  \frac{9q-6 }{10q-12}  \in (1,7/6).\enn

 Indeed, applying Lemma \ref{ADN} to \eqref{ell1} yields   that
 \be\label{ma3.5} \ba  \|\na^2 u\|_{L^q} \le &    C \|\n u_t\|_{L^q}+C\|\n u\cdot \na u\|_{L^q} + C\|\na P\|_{L^q}\\
 \le & C \|\n u_t\|_{L^2}^{\frac{6-q}{2q}} \|\n u_t\|_{L^6}^{\frac{3q-6}{2q}}+C\| \n\|_{L^\infty}  \| u\|_{L^\infty}  \|\na u\|_{L^q}+CM_P(\psi)\psi^\alpha \\
\le & C\psi^{\alpha}\|\sqrt{\n} u_t\|_{L^2}^{\frac{6-q}{2q}}\|\na u_t\|_{L^2}^{\frac{3q-6}{2q}}+C\psi^\alpha  \|\na   u\|_{H^1}^{\frac32}  + CM_P(\psi)\psi^\alpha\\
\le & C\psi^{\alpha}\|\sqrt{\n}u_t\|_{L^2}^{\frac{6-q}{2q}}\|\na u_t\|_{L^2}^{\frac{3q-6}{2q}} +C\psi^\alpha\|\sqrt{\n}u_t\|_{L^2}^{\frac32}+ CM_P^{\frac32}(\psi)\psi^\alpha,
\ea \ee
where in the last inequality one has used \eqref{ell2}.  Combining this with \eqref{ell2},  \eqref{3.1}, and \eqref{li-1a} shows that
 \bnn\label{ma3.9} \ba &\int_0^t  \|\na^2 u\|_{ L^2\cap L^q}^{p_0} ds \\
& \le  C\int_0^t \psi^{\alpha} s^{-p_0/2} \xl(s\|\n^{1/2}u_t\|_{L^2}^2\xr)^{\frac{6-q}{4q}p_0}\xl(s\|\na u_t\|_{L^2}^2\xr)^{\frac{3q-6}{4q}p_0}ds\\
& \quad +C\int_0^t \|\n^{1/2}u_t\|_{L^2}^2ds+ C\int_0^t M_P^{3/2}(\psi)\psi^\alpha ds \\
 & \le C   \exp\left\{C \int_0^tM_P^2(\psi)\psi^\alpha ds\right\} \int_0^t \xl(\psi^{\alpha}+ s^{-\frac{31q^2+12q-36}{26q^2+48q-72}}+s\|\na u_t\|_{L^2}^2\xr) ds\\&\quad+C\exp\left\{C \int_0^tM_P^2(\psi)\psi^\alpha ds\right\}   \\
&  \le C   \exp\left\{C \int_0^tM_P^2(\psi)\psi^\alpha ds\right\} ,
\ea \enn which proves \eqref{bb6-1} and   finishes
 the   proof of    Lemma \ref{l3.4}.   \hfill $\Box$

   Now, we are in a position to prove
  Proposition \ref{pro}.

\emph{Proof of Proposition \ref{pro}}. It follows from  \eqref{3.1}  and \eqref{b1a} that
\bnn\ba \psi(t)
&\le    C_1\exp\left\{C_2\int_{0}^{t}M_P^2(\psi) \psi^{\alpha} ds \right\}.\ea\enn Since $\psi(0)<\ti M  \triangleq  {C_1e} ,$ standard arguments  yield  that for    $T_0\triangleq \min\{1,[C_2M_P^2(\ti M)\ti M^\alpha]^{-1}\} ,$
\be\la{mo2} \sup\limits_{0\le t\le T_0}\psi(t)\le \ti  M,\ee which together with  \eqref{ell2} and  \eqref{li-1a} gives  \be\la{mo1}\ba &\sup\limits_{0\le t\le T_0}t\xl( \|\na^2u\|_{L^2}^{2}+ \|\sqrt{\n} u_t\|_{L^2}^{2}\xr)+  \int_0^{T_0}\left( t \|\na u_t\|_{L^2}^{2} + \|\na^2u\|_{L^2}^{2}\right) dt\le C.\ea\ee

 Next, multiplying \eqref{zb1} by $u_{tt}+u\cdot \na u_t$ and integrating the resulting equation by parts lead to
\be\la{sp9} \ba
& \frac{1}{2}\frac{d}{dt}\int \left(\mu|\nabla u_t|^2 + (\lambda +
\mu)({\rm div}u_t)^2\right)dx+\int_{ }\rho |u_{tt}+u\cdot\na u_t|^2dx
\\
&=\frac{d}{dt}\left(- \int_{
}\rho_t u\cdot\nabla u\cdot u_tdx-\frac{1}{2}\int_{ }\rho_t |u_t|^2 dx+ \int_{ }P_t {\rm
div}u_tdx\right)\\&\quad +
\int_{ } \rho_{tt} u\cdot\nabla u  \cdot u_tdx+
\int_{ }\rho_{t} (u\cdot\nabla u )_t\cdot u_tdx\\&\quad+ \frac{1}{2}\int (\rho_{tt}+\div(u\n_t)) |u_t|^2 dx-
\int_{ } \rho_{t} u\cdot\nabla u \cdot (u\cdot\na  u_t)dx\\ &\quad-\int_{ }\rho
u_t\cdot\nabla u\cdot( u_{tt}+u\cdot\na u_t)dx- \mu\int \pa_iu_t \pa_iu\cdot\na u_tdx\\ &\quad+\frac{\mu}{2} \int\div u|\na u_t|^2dx-(\mu+\lm)\int \div u_t\na u\cdot\na u_tdx\\&\quad+\frac{\mu+\lm}{2}\int \div u(\div u_t)^2dx - \int_{ }P_{tt}{\rm div}u_tdx\\&\quad +\int_{ }P_{t }{\rm div}(u\cdot\na u_t)dx  \triangleq
\frac{d}{dt}I_0+ \sum\limits_{i=1}^{11}I_i. \ea \ee

  We estimate each $I_i (i=0,\cdots,11)$ as follows:

First, it follows from  \eqref{n1}$_1,$   \eqref{mo2},  and \eqref{ell2}   that  \be \ba \la{sp10}|I_0|& =\left|-\frac{1}{2}\int_{
}\rho_t |u_t|^2 dx- \int_{ }\rho_t u\cdot\nabla u\cdot u_tdx+
\int_{ }P_t {\rm div}u_tdx\right|\\ &\le C\left|\int_{ } {\rm
div}(\n
u)|u_t|^2dx\right|+C\norm[L^2]{\rho_t}\|u\|_{L^6}\|\na u\|_{L^6}
\norm[L^6]{u_t}\\&\quad+C\|P_t\|_{L^2}\|\nabla u_t\|_{L^2}\\ &\le C \int_{
} \n |u||u_t||\nabla u_t| dx +C(1+\|\na u\|_{H^1}^2)\|\nabla u_t\|_{L^2} \\
&\le C
\|u\|_{L^6}\|\n^{1/2} u_t\|_{L^2}^{1/2} \|\nabla
u_t\|_{L^2}^{3/2} +C(1+\|\na u\|_{H^1} )\|\nabla u_t\|_{L^2}\\ &\le  \ve\|\nabla
u_t\|_{L^2}^2+C(\ve)\|\n^{1/2} u_t\|_{L^2}^{2}+C,\ea\ee
where in the third inequality we have used \be \la{pn2q}\|\n_t\|_{L^2}+\|P_t\|_{L^2}\le C\|u\|_{L^6}(\|\na \n\|_{L^3}+\|\na P\|_{L^3})+C\|\na u\|_{L^2}\le C.\ee

Next, using  \eqref{n1}$_1$  and  \eqref{mo2},  we have \be\la{pn1q} \|\n_t\|_{L^2\cap L^q}+ \|P_t\|_{L^2\cap L^q}\le C\|\na u\|_{H^1},\ee
which together with \eqref{n1}$_1$  and \eqref{mo2} yields that \be \la{sp11}\ba
  |I_1|&=\left|\int  \rho_{tt} u\cdot\nabla u  \cdot u_{t}dx
 \right|\\
& = \left|  \int_{ } ( \rho_{t } u+\n u_t)\cdot\nabla (u\cdot\nabla u  \cdot u_{t}   )dx\right|\\ &\le C \|  \rho_{t } u+\n u_t \|_{L^3} (\|\na(u\cdot\na u)\|_{L^2}\|u_t\|_{L^6}+\| u\cdot\na u \|_{L^6}\|\na u_t\|_{L^2})\\&\le C\left(\|\na u\|_{H^1}^2+\|\n^{1/2}u_t\|_{L^2}^{1/2}\|\na u_t\|_{L^2}^{1/2}\right)\|\na u\|_{H^1}^2\|\na u_t\|_{L^2}\\&\le C\|\na u\|_{H^1}^2\|\na u_t\|_{L^2}^2+C\|\na u\|_{H^1}^6+C\|\n^{1/2} u_t\|_{L^2}^2\|\na u \|_{H^1}^2 ,\ea \ee
and that
\be \la{sp12}\ba
  |I_2|&=\left|\int_{ }\rho_t\left( u\cdot\nabla u \right)_t\cdot u_{t}dx
 \right|\\
& \le C\|\n_t\|_{L^3}\|(u\cdot\na u)_t\|_{L^2}\|u_t\|_{L^6}\\
& \le C\|\na u\|_{H^1}^2 \|\na u_t\|_{L^2}^2. \ea \ee
Since \eqref{n1}$_1$ implies $\n_{tt}+\div (u\n_t)=-\div(\n u_t),$ we have
 \be \ba
 |I_3|&=
 \frac{1}{2}\left|\int_{ }  \rho u_t \cdot\nabla |u_t|^2 dx\right|\\
 & \le  C
   \norm[L^2]{\rho^{{1/2}}u_t}^{1/2}\|u_t\|_{L^6}^{1/2}
  \|u_t\|_{L^6}
\norm[L^2]{\nabla u_t}\\& \le  C
   \norm[L^2]{\rho^{{1/2}}u_t}^{1/2}
  \|\na u_t\|^{5/2}_{L^2}\\&\le C\|\na u_t\|_{L^2}^2\left( t\|\na u_t\|_{L^2}^2+
   \norm[L^2]{\rho^{{1/2}}u_t}^2+t^{-1/2}\right).
  \ea \ee

Next, Holder's inequality gives
\be  \ba
  |I_4|&=\left|\int_{ } \rho_{t} u\cdot\nabla u \cdot (u\cdot\na  u_t)dx
 \right|\\
& \le C\|\n_t\|_{L^3}\||u|^2|\na u|\|_{L^6}\|\na u_t\|_{L^2}\\
& \le C\|\na u\|_{H^1}^2\|\na u_t\|_{L^2}^2+C\|\na u\|_{H^1}^6, \ea \ee
\be\ba\la{sp13}  |I_5|&= \left| \int_{ }\rho u_t\cdot\nabla
u\cdot (u_{tt}+u\cdot\na u_t) dx\right|  \\& \le   C\|\n^{1/2}(u_{tt}+u\cdot\na u_t)\|_{L^2} \|\n^{1/2}u_t\|_{L^3}\|\na u\|_{L^6}\\& \le   \frac12\|\n^{1/2}(u_{tt}+u\cdot\na u_t)\|_{L^2}^2+C \|\n^{1/2}u_t\|_{L^2}\|\na u_t\|_{L^2}\|\na u\|_{H^1}^2, \ea\ee
and
\be\ba \sum_{i=6}^9|I_i|\le C\|\na u_t\|_{L^2}^2\|\na u\|_{L^\infty}.\ea\ee

Finally, direct calculations together with \eqref{pn1q} lead to
\be \la{sp15}\ba & |I_{10}+I_{11}|\\&=\left| \int_{ }P_{tt}{\rm div}u_tdx-\int_{ }P_{t }{\rm div}(u\cdot\na u_t)dx\right|\\&=\left| \int_{ }P_{tt}{\rm div}u_tdx-\int_{ }P_{t } u\cdot\na{\rm div} u_t dx-\int_{ }P_{t } \na u\cdot\na u_t dx\right|\\&=\left| \int_{ }(P_{tt} + u\cdot\na P_{t } ){\rm div} u_t dx+ \int_{ }P_{t } \div u {\rm div} u_t dx - \int_{ }P_{t } \na u\cdot\na u_t dx\right|\\&
\le C\int\left(|P_t||\na u||\na u_t|+|\na u_t|^2+|u_t||\na P||\na u_t|\right)dx\\&\le C(\|P_t\|_{L^3}\|\na u\|_{H^1}+\|\na P \|_{L^3}\|  u_t\|_{L^6})\|\na u_t\|_{L^2}+C\|\na u_t\|_{L^2}^2\\&\le C \|\na u\|_{H^1}^2 \|\na u_t\|_{L^2}+C\|\na u_t\|_{L^2}^2,\ea\ee
where in the fourth inequality, we have used \be\la{s4} P_{tt}+ u\cdot\nabla P_t=-( \gamma P_t{\rm div}u +
\gamma P{\rm div}u_t + u_t\cdot\nabla P ), \ee due to \eqref{p0}.

Putting all the estimates (\ref{sp11})--(\ref{sp15}) into (\ref{sp9}) and choosing $\ve$ suitably small give
\be\la{nsp19} \ba
  &\Psi' (t) +\int_{ }\rho |u_{tt}+u\cdot\na u_t|^2dx
 \\&\le  C\|\na u_t\|_{L^2}^2\left( t\|\na u_t\|_{L^2}^2+
   \norm[L^2]{\sqrt{\rho} u_t}^2+\|\na u\|_{L^\infty}+ \|\na u\|_{H^1}^2+t^{-1/2}\right)\\&\quad+C\|\na u\|_{H^1}^6 +C\|\sqrt{\n}u_t\|_{L^2}^2\|\na u\|_{H^1}^2+C, \ea \ee
where
\be\ba \label{ma4.2}\notag\Psi (t)\triangleq\mu\|\nabla u_t\|_{L^2}^2+(\mu+\lambda)\|\div u_t\|_{L^2}^2-2I_0
\ea\ee
satisfies
\be\ba \label{ma4.3}    \frac{\mu}{2} \|\na u_t\|_{L^2}^2 -C \|\sqrt{\n}  u_t\|_{L^2}^{2}-C\le \Psi(t)\le C\|\na u_t\|_{L^2}^2 +C  \|\sqrt{\n} u_t\|_{L^2}^{2}+C, \ea\ee
  owing to (\ref{sp10}).  Hence,  multiplying \eqref{nsp19} by $t^2$, we obtain after  using Gronwall's
inequality, \eqref{ma4.3}, (\ref{mo2}),  and (\ref{mo1})   that
\be\ba \label{ma4.1}\sup\limits_{0\le t\le T_0}t^2 \|\na u_t\|_{L^{2}}^2  +\int_{0}^{T_0}t^2 \|\n^{1/2}  u_{tt}\|_{L^{2}}^2 dt \le C,\ea\ee where we have used  the following simple fact that \be\la{md12} \int \n |u|^2|\na u_t|^2dx\le C\|\na u\|_{H^1}^2\|\na u_t\|_{L^2}^2. \ee Combining  \eqref{ma4.1}, \eqref{mo2}, \eqref{mo1}, and  \eqref{ma3.5} gives  \eqref{o1} and completes
the proof of Proposition \ref{pro}.   \hfill $\Box$

\begin{corollary} \la{cor21}Assume that $(\n_0,u_0)$ satisfies \eqref{co2} with some $g\in L^2.$ Then there exists some positive constant  $\ti C$ depending only on       $\mu$, $\lam, $ $P$, $q,$    $\tn,$   $ \psi(0),$ $\|\na u_0\|_{H^1},$ $\|g\|_{L^2},$ and $\O$ if $\O_R=\O $  such that \be\la{m3.12}\ba &\sup_{0\le t\le T_0}\xl(    \|\na u\|_{H^1} +\|\sqrt{\n} u_t\|_{L^2} +t (\|\na u_t\|_{L^2}^{2}+\|\na^2u\|_{L^q}^{2})\xr)\\&+\int_0^{T_0}\|\na u_t\|_{L^2}^2dt\le \ti C. \ea\ee \end{corollary}

 \emph{Proof}.  Taking into account on the compatibility conditions \eqref{co2}, we can define
  \bnn \ba \n^{1/2}u_t(x, t=0)= -g-\n_0^{1/2}u_0\cdot\na u_0,\ea\enn which together with \eqref{a4.6}, \eqref{o1}, and Gronwall's inequality yields \be\la{md3.12} \sup_{0\le t\le T_0}\int\n |u_t|^2dx+\int_0^{T_0}\|\na u_t\|_{L^2}^2dt\le \ti C. \ee
It thus follows from this, \eqref{ell2},   and \eqref{o1} that
 \be\label{mf41} \ba &\sup\limits_{0\le t\le T_0}  \|\na u\|_{H^1} \le   \ti C.\ea\ee which combined with \eqref{nsp19}, \eqref{ma4.3},   \eqref{md3.12}, and \eqref{md12} gives
\be\la{vcs1}\ba \sup\limits_{0\le t\le T_0}t  \|\na u_t\|_{L^{2}}^2  +\int_{0}^{T_0}t \|\n^{1/2}  u_{tt}\|_{L^{2}}^2 dt \le \ti C.\ea\ee  Combining this, \eqref{md3.12}, \eqref{mf41}, and \eqref{ma3.5} gives \eqref{m3.12}  and completes the proof of Corollary \ref{cor21}.  \hfill $\Box$

\section{A priori estimates (II)}\la{sec4}

This section will show some higher order estimates of the solutions with the  initial data satisfying additional compatibility conditions \eqref{co2} and further regularity assumptions  \eqref{1.c1}.
In this section,  the generic positive
constant $C $   depends only   on  $\mu$,  $\lam,$ $P$, $q,$  $\tn,$  $\|\na u_0\|_{H^1}$, and $\|\n_0-\tn\|_{L^{\ti p}\cap D^{1}\cap W^{1,q}}$,  $\|\na^2  \n_0\|_{L^2\cap L^q}, $   $\|\na^2  P(\n_0)\|_{L^2\cap L^q}, $    and $\| g\|_{L^2}.$

\begin{lemma}\label{lem4.5}  It holds that
\begin{equation}\la{5.13a}\ba
 \sup_{0\leq t\leq T_0}\left(  \|\na \n\|_{H^1}+\|\na  P \|_{H^1}+\| \n_t\|_{H^1}+\| P_t\|_{H^1}+t\|\na u\|_{H^2}^2\right) \le C.\ea
\end{equation}
\end{lemma}

\pf It follows from \eqref{n1}$_1,$   \eqref{p0}, and \eqref{o1}    that\be\la{ua2}
 \ba\lefteqn{
\frac{d}{dt}\left(\norm[L^2]{\nabla^2P }  +\norm[L^2]{\nabla^2 \rho
} \right)}\\& & \le C(1+\norm[L^{\infty}]{\nabla
u})\left(\norm[L^2]{\nabla^2P }  +\norm[L^2]{\nabla^2 \rho
} \right) + C\|\na^2u\|_{H^1} . \ea\ee
 Applying Lemma \ref{ADN} to \eqref{ell1} shows
\be\la{sp20} \ba\|\nabla^2 u\|_{H^1}  &\le
 C (\|\n (  u_t+u\cdot\na u)\|_{H^1}+ \|\nabla P\|_{H^1}) \\ &\le C+C \|\na  u_t\|_{L^2}+C \|\na^2 P\|_{L^2},\ea\ee where in the second inequality we have used \eqref{o1},  \eqref{ell2},  and the following simple fact:
  \be \label{jia1} \ba \|\na(\n (u_t+u\cdot\na u))\|_{L^2}
  &\le
 \||\nabla \n | |  u_t|  \|_{L^2}+ \|\n \nabla   u_t  \|_{L^2}+ \|\n|\nabla  u|^2\|_{L^2}
 \\ &\quad
 + \||\nabla \n|| u||\nabla u| \|_{L^2}
 + \|  \n |u || \nabla^2 u| \|_{L^2}\\
 &\le  C
 \|\nabla \n \|_{L^3} \|  u_t  \|_{L^6}+ C\| \nabla   u_t  \|_{L^2}+C\|\na u\|_{H^1}^2
\\
 &\quad  + C\| u\|_{L^\infty}(\| \nabla \n\|_{L^3}\|\nabla u \|_{L^6}
 + C \| \nabla^2 u  \|_{L^2})\\ &\le C+C\| \nabla   u_t  \|_{L^2} \ea\ee due to \eqref{o1} and  (\ref{m3.12}). Using
(\ref{ua2}),   \eqref{sp20}, \eqref{m3.12}, and Gronwall's inequality, one   obtains
\begin{equation}\la{5.1ga}\ba
 \sup_{0\leq t\leq T_0}\left(  \|\na^2\n\|_{L^2}+\|\na^2 P \|_{L^2}+t\|\na^2u\|_{H^1}^2\right) \le C.\ea
\end{equation}

Finally, applying $\na$ to (\ref{p0}) yields \bnn \nabla P_t+u\cdot\nabla\nabla
P+\nabla u\cdot\nabla P+\ga \nabla P {\rm div}u+\ga P  \nabla{\rm
div}u=0,\enn which together with (\ref{5.1ga}),  \eqref{o1},  and (\ref{m3.12}) yields
\be \la{pn3q} \|\nabla P_t\|_{L^2}\le C\|u\|_{L^\infty}\|\nabla^2
P\|_{L^2}+C\|\nabla u\|_{L^6}\|\nabla P\|_{L^3}+C\|\nabla^2
u\|_{L^2}\le C.\ee Similarly, one has $$\|\nabla \n_t\|_{L^2}\le C.$$ Combining this with \eqref{o1}, \eqref{pn2q}, \eqref{pn3q},  and    \eqref{5.1ga} gives \eqref{5.13a}
  and   completes   the proof of
Lemma \ref{lem4.5}. \hfill$\Box$

\begin{lemma}\label{lem4.a5}  It holds that
\begin{equation}\la{5.13n}\ba
 \sup_{0\leq t\leq T_0}\left(\|\nabla^2 \n\|_{L^q } +\|\nabla^2 P  \|_{L^q }\right)  \leq C .\ea
\end{equation}
\end{lemma}

\pf First,
similar to \eqref{ua2}, one has
\be\la{nls21}\ba & (\|\na^2 \n\|_{L^q}+\|\na^2 P\|_{L^q})_t\\&\le  C (1+\|\na u\|_{L^\infty} ) (\|\na^2 \n\|_{L^q}+\|\na^2 P\|_{L^q})   +C  \|\na^2 u\|_{W^{1,q}}    . \ea\ee
Applying Lemma \ref{ADN} to \eqref{ell1} gives
\be \la{sp38}\ba \|\na^2 u\|_{W^{1,q}}& \le C\| \n
 ( u_t+u\cdot \na u) \|_{W^{1,q}}+C\|\na  P
 \|_{W^{1,q}}\\
 &\le C \| \n  ( u_t+u\cdot \na u) \|_{L^2}+C \| \na(\n
( u_t+u\cdot \na u) )\|_{L^{q}} \\&\quad+C\|\na  P
 \|_{L^2} +C\|\na^2  P  \|_{L^{q}}\\  & \le C +C\|\na^2  P
 \|_{L^{q}}  +C \| \na(\n
( u_t+u\cdot \na u) )\|_{L^{q}},\ea\ee due to \eqref{ell2},   (\ref{o1}), and (\ref{m3.12}). For the last term of \eqref{sp38}, it follows from   the Gagliardo-Nirenberg inequality,   \eqref{o1},  \eqref{m3.12},  \eqref{ma3.5},  \eqref{5.13a}, and  \eqref{sp20} that
\be\la{n4.49}\ba   &  \| \na(\n
( u_t+u\cdot \na u) )\|_{L^{q}}\\
&\le C \|\na \n\|_{L^{6q/(6-q)}}(\|  u_t\|_{L^6}+\|  u \|_{L^\infty}\|  \na u \|_{L^6})+C\|\na( u_t+u\cdot \na u)\|_{L^q}
 \\
&\le C(1+\|\na^2\n\|_{L^q})(1+ \|\na  u_t\|_{L^2}) +C\|\na u_t \|_{L^q} \\ & \quad
+C\|\na u \|_{H^1} \| \na  u \|_{H^2}
+C\| u \|_{L^\infty}\|\na^2 u \|_{L^q}\\
&\le C(1+\|\na^2\n\|_{L^q})(1+ \|\na  u_t\|_{L^2}) +C\|\na u_t \|_{L^q} .
    \ea \ee

Then, applying Lemma \ref{ADN} to  \eqref{zb1} yields
  \be\la{sp23}\ba\|\na^2u_t\|_{L^2}
 &\le C\|\n  u_{tt}+\n_t u_t+\n_t u\cdot\nabla u
  +\n u_t\cdot\nabla u+\n u\cdot\nabla u_t+
  \nabla P_t\|_{L^2}   \\
&\le C\left(\|\n  u_{tt}\|_{L^2}+ \|\n_t\|_{L^3}
\|u_t\|_{L^6}+\|\n_t\|_{L^3}\| u\|_{L^\infty}\|\nabla
u\|_{L^6}\right)\\&\quad
  +C\left(\| u_t\|_{L^6}\|\nabla u\|_{L^3}+ \| u\|_{L^\infty}
  \|\nabla u_t\|_{L^2}+\|
  \nabla P_t\|_{L^2}\right)\\ &\le C\|\n^{1/2}  u_{tt}\|_{L^2} +C\|\nabla  u_t\|_{L^2}+C,\ea \ee where in the last inequality we have used  (\ref{m3.12}),  \eqref{o1},    (\ref{pn1q}), and (\ref{5.13a}).
  Combining this with \eqref{m3.12} and \eqref{vcs1} shows
  \be\label{ma4.5}\ba \int_0^{T_0}\|\na u_t \|_{L^q}dt&\le  C\int_0^{T_0}\|\na u_t \|_{L^2}^{(6-q)/(2q)}
\|\na u_t \|_{H^1}^{3(q-2)/(2q)}dt
\\&\le C+C  \int_0^{T_0}t^{-1/2}(t
\|\n^{1/2}  u_{tt} \|^2_{L^2})^{3(q-2)/(4q)}dt\\&\le C+C  \int_0^{T_0}\left(t^{-2q/(q+6)}+ t
\|\n^{1/2}  u_{tt} \|^2_{L^2}\right)  dt\le C.\ea\ee

  Finally, putting \eqref{sp38} and  \eqref{n4.49} into \eqref{nls21}  and using Gronwall's inequality, \eqref{m3.12}, and \eqref{ma4.5},  we  obtain   \eqref{5.13n} and complete the proof of Lemma \ref{lem4.a5}.  \hfill$\Box$

\begin{lemma} \la{4.m44} It holds that
\begin{equation}\label{4.m30}\ba
&\sup_{0\leq t\leq T_0} t\left(  \|\na^3
u \|_{  L^q} + \|\na
u_t \|_{H^1} +\|\sqrt{\rho }
u_{tt}\|_{L^2} \right)  +\int_0^{T_0}  t^2  \|\nabla u_{tt}\|_{L^2}^2  dt\leq
C .\ea
\end{equation}
\end{lemma}

\pf  We claim that
\begin{equation}\la{5.o2}
\sup_{0\leq t\leq T_0}t^2 \|\sqrt{\rho }
u_{tt}\|_{L^2}^2+\int_0^{T_0}  t^2 \|\nabla u_{tt}\|_{L^2}^2 dt  \leq
C ,
\end{equation}
which together with    \eqref{m3.12} and \eqref{sp23} yields  that
\begin{equation}\la{5.p6}
\sup_{0\leq t\leq T_0}t \|\na   u_t\|_{H^1}  \leq
C .
\end{equation}
It thus follows from this,  \eqref{sp38},    \eqref{n4.49}, and \eqref{5.13n}  that
\begin{equation}\la{5.p5}
\sup_{0\leq t\leq T_0}t\|\na^3  u\|_{L^q}\leq
C.
\end{equation}    Combining \eqref{5.o2}--\eqref{5.p5} yields  \eqref{4.m30}.

Now, it remains to prove \eqref{5.o2}. In fact, differentiating   (\ref{zb1})   with
respect to $t$ leads to
 \be\la{sp30}\ba &\n u_{ttt}+\n u\cdot\na u_{tt}-\mu\Delta
u_{tt}-(\mu+\lambda)\nabla{\rm div}u_{tt}\\&= 2{\rm div}(\n u)u_{tt}
+{\rm div}(\n u)_{t}u_t-2(\n u)_t\cdot\na u_t-(\n_{tt} u+2\n_t u_t)
\cdot\na u\\& \quad- \n u_{tt}\cdot\na u-\na P_{tt}.
 \ea\ee
Multiplying \eqref{sp30} by $u_{tt}$  and integrating the resulting equation   by parts   yield
\be \la{sp31}\ba &\frac{1}{2}\frac{d}{dt}\int_{ }\n
|u_{tt}|^2dx+\int_{ }\left(\mu|\na u_{tt}|^2+(\mu+\lambda)({\rm
div}u_{tt})^2\right)dx \\&=-4\int_{ }  \n u\cdot\na
 u_{tt}\cdot u_{tt} dx-\int_{ }(\n u)_t\cdot \left[\na (u_t\cdot u_{tt})+2\na
u_t\cdot u_{tt}\right]dx\\&\quad -\int_{
}(\n_{tt}u+2\n_tu_t)\cdot\na u\cdot u_{tt}dx-\int_{ }   \n
u_{tt}\cdot\na u\cdot  u_{tt} dx\\&\quad+\int_{ } P_{tt}{\rm
div}u_{tt}dx\triangleq\sum_{i=1}^5K_5.\ea\ee

Using  \eqref{o1},  (\ref{m3.12}), and  (\ref{5.13a}), we can estimate each $K_i(i=1,\cdots,5)$ as follows: \be \la{sp32} \ba |K_1|&\le
C\|\n^{1/2}u_{tt}\|_{L^2}\|\na u_{tt}\|_{L^2}\| u \|_{L^\infty}\\
&\le \ve \|\na u_{tt}\|_{L^2}^2+C(\ve)\|\n^{1/2}u_{tt}\|^2_{L^2},\ea\ee
 \be \la{sp33}\ba |K_2|&\le C\left(\|\n
u_t\|_{L^3}+\|\n_t u\|_{L^3}\right)\left(\| u_{tt}\|_{L^6}\| \na
u_t\|_{L^2}+\| \na u_{tt}\|_{L^2}\| u_t\|_{L^6}\right)\\
&\le
C\left(\|\n^{1/2} u_t\|^{1/2}_{L^2}\|u_t\|^{1/2}_{L^6}+\|\n_t
\|_{L^6}\| u\|_{L^6}\right)  \| \na u_{tt}\|_{L^2}  \| \na u_{t}\|_{L^2} \\
&\le
C\left( \|\na u_t\|_{L^2}+1\right)  \| \na u_{tt}\|_{L^2}  \| \na u_{t}\|_{L^2} \\
&\le \ve
\|\na u_{tt}\|_{L^2}^2+C(\ve)\| \na u_{t}\|_{L^2}^4+C(\ve),\ea\ee

\be  \la{sp34}\ba |K_3|&\le C\left(\|\n_{tt}\|_{L^2}
\|u\|_{L^\infty}\|\na u\|_{L^3}+\|\n_{
t}\|_{L^6}\|u_{t}\|_{L^6}\|\na u \|_{L^2}\right)\|u_{tt}\|_{L^6} \\
&\le \ve \|\na u_{tt}\|_{L^2}^2+C(\ve)\|\n_{tt}\|_{L^2}^2+C(\ve)\| \na u_{t}\|_{L^2}^2,\ea\ee and
\be  \la{sp36}\ba |K_4|+|K_5|&\le C\|\n u_{tt}\|_{L^2} \|\na
u\|_{L^3}\|u_{tt}\|_{L^6} +C \|P_{tt}\|_{L^2}\|\na
u_{tt}\|_{L^2}\\
&\le \ve \|\na u_{tt}\|_{L^2}^2+C(\ve)\|\n^{1/2}u_{tt}\|^2_{L^2}
+C(\ve)\|P_{tt}\|^2_{L^2}. \ea\ee
Substituting (\ref{sp32})--(\ref{sp36}) into (\ref{sp31})
and choosing $\ve$ suitably small lead to
 \be\la{nsp37}\ba &\frac{d}{dt}\|\n^{1/2}u_{tt}\|^2_{L^2}+\mu\|\na u_{tt}\|_{L^2}^2\\&\le C \|\n^{1/2}u_{tt}\|^2_{L^2}+ C \| \na u_{t}\|_{L^2}^4+C+C\|\n_{tt}\|^2_{L^2}+C\|P_{tt}\|^2_{L^2}.\ea\ee

 Finally, it follows from \eqref{s4}, (\ref{5.13a}),      and (\ref{md3.12})  that \be \la{ntt1} \ba  \int_0^{T_0} \|P_{tt}\|_{L^2}^2ds  & \le C\int_0^{T_0}
\left( \|u\|_{L^\infty}\|\nabla
P_t\|_{L^2}+\|P_t\|_{L^6}\|\nabla u\|_{L^3}\right)^2dx\\&\quad+C\int_0^{T_0}
\left( \|\nabla
u_t\|_{L^2} +\|u_t\|_{L^6}\|\nabla P\|_{L^3} \right)^2dt \le C.\ea\ee
Similarly, one has \be \la{ntt2}\int_0^{T_0}\|\n_{tt}\|_{L^2}^2dt\le C.\ee
Multiplying \eqref{nsp37} by $t^2$ and using \eqref{m3.12}, \eqref{vcs1},
  (\ref{ntt1}), and
(\ref{ntt2}),  we obtain  \eqref{5.o2} and     finish  the proof of Lemma \ref{4.m44}.  \hfill$\Box$

\section{Proofs of Theorems \ref{t1} and \ref{t3} }

To prove Theorems \ref{t1}--\ref{t3},
we will only deal with the case that $\O$ is  bounded.  Since for the Cauchy problem, all the a priori estimates obtained in sections 3 and 4 are independent of the radius $R$, one can use the standard domain expansion technique to treat the whole space case, please refer to  \cite{L1} and references therein.

{\it Proof of Theorem  \ref{t1}.}
Let  $(\n_{0},u_{0})$  be as  in Theorem \ref{t1}. For $\de>0,$  we choose $0\le\hat\n_0^\de \in C^\infty(\O)  $ and $u_0^\de\in C^\infty_0(\O)$ satisfying
      \be\label{ma5.3}  \lim_{\de\rightarrow 0}\left(\|\hat\n_0^\de-\n_0\|_{W^{1,q}}+\|u_0^\de-u_0\|_{H^1}\right) =0. \ee
Then, in terms of Lemma \ref{th0}, the problem  \eqref{n1}--\eqref{n5} with the initial data $(\hat\n_0^{\delta}+\de,(\hat\n_0^{\delta}+\de)u_0^{\delta})$ has  a unique smooth solution   $(\n^{\delta},u^{\delta}) $ on $\O\times [0,T_{\de}] $ for some $T_\de>0.$  Moreover,  Proposition \ref{pro} shows that there exist   two positive constants  $T_0$ and $M$ independent of   $\delta$ such that \eqref{o1} holds for $(\n^{\delta},u^{\delta}).$  More precisely, it holds

  \be\label{ma5.5}\ba&\sup\limits_{0\le t\le T_0}\xl(\|\na u\|_{L^2} +    \|  \n \|_{ H^{1}\cap W^{1,q}}+ \| P(\n)  \|_{  H^{1}\cap W^{1,q}}+t(\|\na^2u\|_{L^2}^{2}+\|\sqrt{\n} u_t\|_{L^2}^{2})\xr)\\&+\sup\limits_{0\le t\le T_0}\xl(t^2(\|\na u_t\|_{L^2}^{2}+\|\na^2u\|_{L^q}^{2})\xr) +\int_0^{T_0}t\|\na u_t\|_{L^2}^2dt \le M\ea\ee
 \be \ba &\sup\limits_{0\le t\le T_0}\xl(\|\n^{\delta}\|_{W^{1,q}}+ \|\n_t^{\delta}\|_{L^2}+\|u^{\delta}\|_{H^1}\ +t^{1/2}\|\na^2u^{\delta}\|_{L^2}+\|\n^{\delta}u^{\delta}\|_{H^1}\xr) \\
 &\quad+  \int_0^{T_0} \left(\norm[L^{q}]{\nabla^2u^{\delta}}^{p_0}+t\|\na u^{\delta}_t\|_{L^2}^{2}+t\|\na^2u^{\delta}\|_{L^q}^{2}+ \|\na^2u\|_{L^2}^{2}+\|(\n^{\delta}u^{\delta})_t\|_{L^2}^2\right)dt\le \bar C, \ea\ee
 where $\bar C$ is independent of   $\delta$.
  With all the estimate \eqref{ma5.5}  at hand,  we  find  that the sequence
$(\n^{\delta},u^{\delta})$ converges, up to the extraction of subsequences, to some limit $(\n,u)$   in the
obvious weak sense. That is, as   $\delta\rightarrow 0$, we have
\be\la{kq3} \n^{\delta} \rightarrow  \n,\mbox{in } L^\infty(0,T_0;  L^\infty),\ee
\be  \n^{\delta} \rightharpoonup  \n,\mbox{ weakly * in } L^\infty(0,T_0;  W^{1,q}),\ee
\be u^{\delta} \rightharpoonup u,\mbox{ weakly * in } L^\infty(0,T_0;  H^1),\ee
\be  \na^2 u^{\delta}\rightharpoonup \na^2 u ,\mbox{ weakly  in }L^{p_0}(0,T_0; L^q)\cap L^2(\O\times (0,T_0)),\ee
\be  t^{1/2}\na^2 u^{\delta}\rightharpoonup t^{1/2}\na^2 u ,\mbox{ weakly  in }L^2(0,T_0; L^q ),\ee
 \be t^{1/2} \na  u^{\delta}_t\rightharpoonup t^{1/2}\na  u_t ,\mbox{ weakly  in } L^2(\O\times (0,T_0)),\ee
 \be \la{kq4}  \n^{\delta}u^{\delta} \rightarrow  \n u,\mbox{in } L^\infty(0,T_0;  L^2). \ee
 Then letting  $\delta\rightarrow 0$, it follows from \eqref{kq3}-\eqref{kq4} that
$(\n,u)$ is a strong solution of \eqref{n1}-\eqref{n5} on $\O\times (0,T_0]$ satisfying \eqref{1.10}.  The  proof of the existence part of Theorem \ref{t1} is finished.

It only remains to prove the uniqueness of the strong solutions  satisfying \eqref{1.10}. Indeed, we will use the method which is due to Germain \cite{Ge}.   Let   $(\rho,u)$ and $(\bar\rho,\bar u)$ be two strong solutions satisfying \eqref{1.10}  with the same initial data. Subtracting the mass   equation  for $(\rho,u)$ and $(\bar\rho,\bar u)$ gives
\be\la{5.2}
 H_t + \bar u\cdot\nabla H +H\div \bar u + \rho \div U  + U\cdot\nabla \rho= 0,
\ee
with
$$ H\triangleq\rho-\bar\rho, \quad U\triangleq u-\bar u.$$
For $3/2\le r\le 2$, multiplying   \eqref{5.2} by $rH |H|^{r-2}$ and  integrating the resulting equation by parts lead to
  \be\label{ma5.6} \ba
  \frac{d}{dt}\|H\|_{L^r}^r&\le C\int \div \bar u |H|^rdx+C\int \n |\na U| |H|^{r-1}dx+C\int |U||\na \n| |H|^{r-1}dx\\
  &\le C\|\na \bar u\|_{L^\infty}\|H\|_{L^r}^r+C\xl(\|\n\|_{L^{\frac{2r}{2-r}}}+\|\na \n\|_{L^{\frac{6r}{6-r}}}\xr)\|\na U\|_{L^2}\|H\|_{L^r}^{r-1}\\
  &\le C\|\na \bar u\|_{L^\infty}\|H\|_{L^r}^r+C\|\na U\|_{L^2}\|H\|_{L^r}^{r-1},
  \ea\ee
where one has used $\n\in H^1\cap W^{1,q}$. This together with Gronwall's inequality and \eqref{o1} gives
  \be\label{ma5.7} \ba
 \|H\|_{L^r} \le C\int_0^t\|\na U\|_{L^2}ds,~~~\mbox{for}~3/2\le r\le 2.
  \ea\ee

   Next, subtracting the momentum equations for $(\rho,u)$ and $(\bar\rho,\bar u)$ yields
\be\ba\la{5.5}
& \rho U_t + \rho u\cdot\nabla U -\mu\lap U - (\mu+\lambda)\na\left(\div  U\right)\\
 &=  - \rho U\cdot\nabla\bar u-H(\bar u_t+\bar u\cdot\nabla\bar u)-\na \left( P(\n)-P(\bar{\n}) \right) ,
\ea\ee
Multiplying (\ref{5.5}) by $U$ and integrating the resulting equations by parts lead to
\be\ba\la{5.6}
& \frac{d}{dt}\int \rho |U|^2dx  + 2\mu \int |\na U|^2 dx\\
 & \le   C \norm[L^{\infty}]{\nabla \bar{u}}\int \n |U|^2dx+C\int|H || U| \left(| \bar u_{t}|+ | \bar u ||\na\bar{u}|\right) dx \\
 &\quad+C\|P(\n)-P(\bar{\n})\|_{L^2}\norm[L^{2}]{\div U}\\
  & \le    C \norm[L^{\infty}]{\nabla \bar{u}}\int \n |U|^2dx+C\|H\|_{L^{3/2}}\|U\|_{L^6}\|\bar u_t\|_{L^6}\\
 &\quad +C\|H\|_{L^2}\|U\|_{L^6}\|\bar u \|_{L^6}\|\na \bar u \|_{L^6}+C\|H\|_{L^2}\|\na U\|_{L^2}\\
    & \le    C \norm[L^{\infty}]{\nabla \bar{u}}\int \n |U|^2dx +C \xl(1+\|\na \bar u_t\|_{L^2}+\|\na^2 \bar u \|_{L^2}\xr)\|\na U\|_{L^2} \int_0^t\|\na U\|_{L^2}ds\\
       & \le    C \norm[L^{\infty}]{\nabla \bar{u}}\int \n |U|^2dx +C \xl(1+t\|\na \bar u_t\|_{L^2}+t\|\na^2 \bar u \|_{L^2}\xr)  \int_0^t\|\na U\|_{L^2}^2ds+\mu \|\na U\|_{L^2}^2\\
     & \le    C\xl(1+t\|\na \bar u_t\|_{L^2}^2 +\|\na \bar u\|_{L^\infty}\xr)  \xl(\int \n |U|^2dx+\int_0^t\|\na U\|_{L^2}^2dt\xr) +\mu \|\na U\|_{L^2}^2
\ea\ee owing to \eqref{o1} and \eqref{ma5.7}. This   together with Gronwall's inequality and \eqref{o1} gives  $U(x,t)=0 $ for  almost everywhere $(x,t)\in \O\times(0,T_0).$  Then, \eqref{ma5.7} implies that $H(x,t)=0 $ for almost everywhere $(x,t)\in \O\times(0,T_0).$ The proof of Theorem  \ref{t1} is completed.  \hfill$\Box$

{\it Proof of Theorem  \ref{t3}.}
Let  $(\n_{0},u_{0})$  be as  in Theorem \ref{t3}, we construct
$\n_0^{\delta}=\hat\n_0^{\delta}+\delta$ where  $0\le \hat\n_0^{\delta}\in  C_0^\infty (\O)$ satisfies \eqref{ma5.3}   and
$$\na^2 \hat\n_0^{\delta}\rightarrow  \na^2\n_{0},~~\na^2 P(\hat\n_0^{\delta})\rightarrow  \na^2 P(\n_{0}),~\quad {\rm in}\,\, L^2\cap L^q,~~~ {\rm as}~ \delta\rightarrow 0.$$
Thus, we have
      \be\label{ma6.1} \ba\begin{cases} \n_0^{\delta}\rightarrow  \n_{0}\quad {\rm in}\,\, W^{1,q}(\O) ,\\ \na^2 \n_0^{\delta}\rightarrow  \na^2\n_{0}\quad {\rm in}\,\, L^2\cap L^q,\\
      \na^2 P(\n_0^{\delta})\rightarrow  \na^2 P(\n_{0})\quad {\rm in}\,\, L^2\cap L^q, \end{cases} ~ {\rm as}~ ~\delta\rightarrow 0.\ea\ee

Then, we consider the unique smooth solution $u_0^{\de}$ of the following elliptic problem:
 \be \la{bb12} \begin{cases}- \mu\lap u^{\delta}_0 - (\mu + \lam)\nabla \div u^{\delta}_0 +  \nabla P(\n_0^{\delta})= \sqrt{\n_0^{\delta}} g^{\delta}, & {\rm in} \,\,  \O,\\ u_{0}^{\delta} =0, &{\rm on} \,\,\partial \O,\end{cases} \ee
where $  g^{\delta}=  g *j_{\de} $  with $j_\de$
being the standard mollifying kernel of width $\de.$

 Subtracting the equations \eqref{co2} and \eqref{bb12} gives
  \be \label{ma6.2}\begin{cases}- \mu\lap \xl(u^{\delta}_0-u_0\xr) - (\mu + \lam)\nabla \div \xl(u^{\delta}_0-u_0\xr) =F, & {\rm in} \,\,  \O,\\  u^{\delta}_0-u_0 =0, &{\rm on} \,\,\partial \O,\end{cases} \ee
  with
  \be \label{ma6.3}\notag F\triangleq-\nabla \xl(P(\n_0^{\delta})-P(\n_0)\xr)+\sqrt{\n_0^{\delta}} g^{\delta}-\sqrt{\n}g.\ee
  Multiplying  \eqref{ma6.2} by $u^{\delta}_0-u_0$, we obtain after integration by parts that
     \be \ba\label{ma6.4}&\|\nabla \xl( u_0^{\delta} - u_0\xr)\|_{L^2}\\
    & \le C\|P(\n_0^{\delta})-P(\n_0)\|_{L^2}+ C\| \sqrt{\n_0^{\delta}} - \sqrt{\n_0} \|_{L^3}+ C\|g^{\delta}- g\|_{L^2}\\
     &\rightarrow 0,~~{\rm as}~ \delta\rightarrow 0,\ea\ee
  due to \eqref{ma5.3} and \eqref{ma6.1}. Moreover,  Lemma \ref{ADN} combined with  \eqref{ma6.2} yields  that
      \be \ba\label{ma6.5}&\|\nabla^2 \xl( u_0^{\delta} - u_0\xr)\|_{L^2}\\
    & \le C\|\na P(\n_0^{\delta})-\na P(\n_0)\|_{L^2}+ C\| \sqrt{\n_0^{\delta}} - \sqrt{\n_0} \|_{L^\infty}+ C\|g^{\delta}- g\|_{L^2}\\
     &\rightarrow 0,~~{\rm as}~  \delta\rightarrow 0,\ea\ee
 owing  to \eqref{ma5.3} and \eqref{ma6.1}.

For the problem  \eqref{n1}-\eqref{n5} with the initial data $(\n_0^{\delta},u_0^{\delta})$ satisfying \eqref{ma5.3} and \eqref{ma6.1}--\eqref{bb12},   Lemma \ref{th0} shows that there exists   a classical solution   $(\n^{\delta},u^{\delta}) $ on $\O\times [0,T_0].$  Moreover,   we deduce from  \eqref{o1} and Lemmas 4.1--4.5 that   the sequence
$(\n^{\delta},u^{\delta})$ converges weakly, up to the extraction of subsequences, to some limit $(\n,u)$   satisfying   \eqref{1.10}, \eqref{ma2}, and \eqref{1.a10}. Moreover, standard arguments yield that $(\n,u)$   in fact is a classical solution to   the problem \eqref{n1}-\eqref{n5}. The proof of Theorem  \ref{t3} is completed.  \hfill$\Box$

\end{document}